\documentclass{amsart}
\usepackage{amsmath}
\usepackage{amssymb}
\usepackage[mathscr]{eucal}
\usepackage{amsthm}
\usepackage{amsfonts}
\usepackage{hyperref}
\usepackage{extarrows}

\usepackage[margin=1in]{geometry}

\makeatletter
\newtheorem*{rep@theorem}{\rep@title}
\newcommand{\newreptheorem}[2]{%
\newenvironment{rep#1}[1]{%
 \def\rep@title{#2 \ref{##1}}%
 \begin{rep@theorem}}%
 {\end{rep@theorem}}}
\makeatother

\newtheorem{lemma}{Lemma}
\newtheorem{theorem}{Theorem}
\newreptheorem{theorem}{Theorem}

\def\bea{\begin{eqnarray}}
\def\eea{\end{eqnarray}}
\def\bean{\begin{eqnarray*}}
\def\eean{\end{eqnarray*}}
\def\be{\begin{equation}}
\def\ee{\end{equation}}

\newcommand{\Gr}{\mathrm{Gr}}

\newcommand{\mZ}{\mathbb{Z}}
\newcommand{\mP}{\mathbb{P}}

\newcommand{\mC}{\mathbb{C}}

\newcommand{\bd}{\mathbf{d}}

\newcommand{\cA}{\mathcal{A}}

\newcommand{\cM}{\mathcal{M}}

\newcommand{\cO}{\mathcal{O}}

\newcommand{\cF}{\mathcal{F}}
\newcommand{\cI}{\mathcal{I}}

\newcommand{\fI}{\mathfrak{I}}

\DeclareMathOperator{\ch}{ch}

\DeclareMathOperator{\Tr}{Tr}

\DeclareMathOperator{\Span}{span}
\DeclareMathOperator{\Aut}{Aut}
\DeclareMathOperator{\Tod}{Tod}
\DeclareMathOperator{\Mor}{Mor}
\DeclareMathOperator{\Def}{Def}

\begin{document}
\title{Poincar\'e polynomials of moduli spaces of stable maps into flag manifolds}
\author{Xiaobo Zhuang}

\date{}

\begin{abstract}
By using Bialynicki-Birula decomposition for the stack of genus zero stable maps to flag manifolds\cite{Oprea2006Tautological}. We calculate the Poincar\'{e} polynomial of the moduli space in degree one and degree two.
\end{abstract}

\keywords{Bialynicki-Birula decomposition, Poincar\'e polynomial}

%\ccode{Mathematics Subject Classification 2000: 58J26, 14K25, 14M15}

\maketitle

\section{Introduction}
In enumerative geometry, when one wants to know about rational curves of degree $\bd\in H_{2}(X)$ in a space $X$, we consider the space $\Mor_{\bd}(\mP^1,X)$ of morphisms from $\mP^1$ to $X$ of degree $\bd\in H_2(X,\mZ)$, and  use intersection theory on $\Mor_{\bd}(\mP^1,X)$ to solve the enumerative problem. However, the problem is that the space $\Mor_{\bd}(\mP^1,X)$ is not compact, so we compactify it. When $X=Fl(r_1,\cdots,r_{l+1};k)=Fl(r_1,\cdots,r_{l+1})$ the flag manifold which parametrizes successive subspaces in $\mC^k$:
  $$
  V_1 \subset V_2 \subset \cdots \subset V_l \subset \mC^k
  $$
with $\dim V_j=\sum_{i=1}^j r_i$ and $\sum_i^{l+1} r_i=k$. There is a natural compactification called the hyperquot scheme, it is wildly used in Gromov-Witten theory and quantum cohomology ring of Grassmannian  and flag manifolds. In \cite{Str1987On} Str{\o}mme derive an implicite formula for the Betti numbers of the Quot schemes using Bialynicki-Birula decomposition for Quot scheme. Later in \cite{Chen2000Poincare}, Chen generalized the method to partial flag manifolds and computed the generating function for the Poincar\'e polynomials of hyperquot schemes.

However there is another natural compactification of $\Mor_{\bd}(\mP^1,X)$, that is Kontsevich's moduli space of stable maps $\cM_0(Fl(r_1,\cdots,r_{l+1}),\bd)$. In \cite{Fulton1996Notes} Fulton and Pandharipande shows that its coarse moduli space is a projective normal variety with an orbifold structure. In \cite{Manin1998Stable}, Manin calculate its virtual Poincar\'e polynomial.  Oprea's work \cite{Oprea2006Tautological} shows that there is a Bialynicki-Birula decomposition for the moduli space of stable maps to projective spaces. Applying Oprea's decomposition to the moduli space of stable maps to Grassmannian. Agrawal  \cite{agrawal2011euler} computes the Euler characteristics of the coarse moduli space of stable maps to Grassmannian in lower degrees and later in \cite{Mart2013Poincar} Mart\'in computes its Poincar\'e polynomial in lower degree. Edwards \cite{edwards2013genus} computes the Euler characteristics of the coarse moduli space of stable maps to flag manifolds in lower degrees.

In this paper, we carry out localization analysis on flag manifolds, and compute the Poincar\'e polynomial of the moduli space of stable maps to flag manifolds in lower degree using the corresponding Bialynicki-Birula decomposition. In our computation, holomorphic Lefschetz formula plays a important role. Our main results are summarized in the following theorem:

\begin{theorem}
  Let $Fl=Fl(r_1,\cdots,r_{l+1})$ be the partial flag manifold, and $\cM_0(Fl(r_1,\cdots,r_{l+1}),\bd)$ be the moduli space of genus zero stable maps of degree $\bd$. Then its Poincar\'e polynomials in degree one and two are:
  \begin{align}
    P_{\cM_0(Fl(r_1,\cdots,r_{l+1}),\check H_i)}(q)&=\frac{[r_i]_{t}[r_{i+1}]_{t}}{1+t}P_{Fl}(q) \\
    P_{\cM_0(Fl,\check H_i+\check H_j)}(q)&=(1+t^2)\frac{[r_i]_{t}[r_{i+1}]_{t}[r_j]_{t}[r_{j+1}]_{t}}{1+t}P_{Fl}(q),\quad j-i>1\\
      P_{\cM_0(Fl,2\check H_i)}(q)
  &=\frac{(1-t^{r_i})(1-t^{r_{i+1}})((1+t^{r_i+r_{i+1}})(1+t^3)-t(1+t)(t^{r_i}+t^{r_{i+1}}))}{(1-t)^2(1-t^2)^2}P_{Fl}(q)
  \end{align}
  where $t=q^2$ and $P_{Fl}(q)=\binom{k}{r_1,\cdots,r_{l+1}}_{q^2}$ is the Poincar\'e polynomial of the partial flag. When the flag is a complete flag, we have:
  \begin{equation}
      P_{\cM_0(Fl,\check H_i+\check H_{i+1})}(q)=\frac{1+2t+3t^2+3t^3+t^4}{(1+t)(1+t+t^2)}P_{Fl}(q)
  \end{equation}
\end{theorem}

We make a commment on the Bialynicki-Birula decomposition used in our paper. In \cite{Oprea2006Tautological}, Oprea consider Bialynicki-Birula decompostion for Deligne-Mumford stacks. Later in \cite{Skowera2013Bia} Skowera extended Oprea's result to that any smooth, proper, tame Deligne-Mumford stack, whose coarse moduli space is a scheme admits a Bialynicki-Birula decomposition. He also remarks that (see Remark 3.6) when the coarse moduli space is a projective scheme then the induced decomposition is filterable in the sense of Oprea, and then one may use Lemma 6 in \cite{Oprea2006Tautological} to compute the the Betti numbers of the moduli space from the fixed locus.

This paper is organized as follows: In Section 2, we recall Bialynicki-Birula decomposition and holomorphic Lefschetz formula. In Section 3, we carry out localization analysis on the moduli space. In section 4, we use holomorphic Lefschetz formula to compute the contributions of fixed locus and prove our main theorems.

\medskip\noindent
{\bf Acknowledgment}: The author thank Professors Jian Zhou and Kefeng Liu for their encouragements, comments and suggestions. He thanks Professor Dragos Oprea for telling him Edward's work \cite{edwards2013genus} and some reference typos. He also wants to thank Xiaowen Hu for helpful discussions.

\section{Preliminaries}
\subsection{Bialynicki-Birula decomposition}
The theory of Bialynicki-Birula decomposition is developed in \cite{Bialynicki1973Some}\cite{Bia1976Some}( see also the book \cite{Carrell2002Torus}). Let $X$ be a smooth projective variety, and $T$ an algebraic torus of dimension one such that $T$ acts on $X$. Suppose the fixed point set $X^T$ is nontrivial. Let $Y_1,\cdots,Y_r$ be the irreducible components of $X^T$. There is a decomposition of the tangent bundle when restricted to $Y_i$:
\begin{equation}
  TX|_{Y_i}= T_i^+\oplus T_i^0 \oplus T_i^-
\end{equation}
where $T_i^+$, $T_i^0$ and $T_i^-$ are subbundles of $TX|_{Y_i}$ such that the torus acts on it with positive, trivial and negative weights respectively. We denote the rank of $T^+_i$ by $p_i$ and that of $T_i^-$ by $n_i$. We define:
\begin{equation}
  Y_i^+=\{x\in X|\lim_{t\to 0} t\cdot x \in Y_i \}
\end{equation}
which are called the {\it plus cells};
\begin{equation}
  Y_i^-=\{x\in X|\lim_{t\to \infty} t\cdot x \in Y_i \}
\end{equation}
which are called the {\it minus cells}. Then we have the {\it plus decomposition}:
\begin{equation}
  X=\cup_{1\leq i\leq r} Y_i^+
\end{equation}
and the {\it minus decomposition}:
\begin{equation}
  X=\cup_{1\leq i\leq r} Y_i^-
\end{equation}
Then (see \cite[Theorem 4.1]{Bialynicki1973Some} or \cite[Theorem 4.2]{Carrell2002Torus}):
\begin{itemize}
  \item each irreducible component $Y_i$ is smooth and the plus(resp. minus) cells are locally closed;
  \item the natural projections morphisms $\pi^+_i: Y^+_i \to Y_i$(resp. $\pi^-_i: Y^-_i \to Y_i$) are $T-$isomorphic to $p_i^+: T_i^+\to Y_i$(resp. $p_i^-: T_i^-\to Y_i$);
  \item $T_i^0\cong TY_i$.
\end{itemize}

And we may use homology "basis" theorem to compute the homology of $X$:
\begin{theorem}[Homology basis theorem]
  \begin{align}
    H_{m}(X)&\cong  \oplus_{i} H_{m-2p_i}(Y_i) \\
    &\cong \oplus_{i} H_{m-2n_i}(Y_i)
  \end{align}
\end{theorem}
So the Poincar\'e polynomial of the total space $X$ can be computed from that of the fixed locus:
\begin{equation}\label{eqn: poincare polynomial of scheme}
  P_X(t)=\sum_i t^{2p_i}P_{Y_i}(t)=\sum_{i} t^{2n_i}P_{Y_i}(t)
\end{equation}
In many cases, the $Y_i$'s are isolated points in $X$. So to compute the Poincar\'e polynomial, it suffices to determine the numbers $p_i$ or $n_i$.

\subsection{Generalization to Deligne-Mumford stacks}
Let $\cM$ be a Deligne-Mumford stack with a one-dimensional torus $T$ acting on it. Let $\cF_i$ be the components of the fixed locus of the action. Let $p_i$ and $n_i$ be the rank of the corresponding subbundles of the tangent bundle as above. In \cite{Oprea2006Tautological},  Oprea proved the Bialynicki-Birula decomposition for Deligne-Mumford stacks provided that  there exists a T-equivariant, affine, \'etale atlas. When the decomposition is filterable, Oprea proves a homology "basis" theorem for Deligne-Mumford stacks:
\begin{theorem}\label{thm: homology basis theorem for DM stack}
  When the decomposition is filterable, the Betti numbers $h^m(\cM)$ of $\cM$ can be computed as
\begin{equation}
  h^i(\cM)=\sum_i h^{i-2n_i}(\cF_i)
\end{equation}
Here $n_i$ is the codimension of $F_i^+$ which equals the number of negative weights on the
tangent bundle of $\cM$ at a fixed point in $\cF_i$.
\end{theorem}
Note that the betti numbers are defined using the cohomology theory in \cite{Behrend2005Cohomology}, and equals the betti number of the coarse moduli space (see \cite[Proposition 36]{Behrend2005Cohomology}). When the Delign-Mumford stack has a projective coarse moduli space, the decomposition is filterable, then the homology basis theorem of Oprea applies (see \cite[Remark 3.6]{Skowera2013Bia}). This fact is used in \cite{agrawal2011euler} and \cite{Mart2013Poincar} to compute the Euler characteristic and Poincar\'e polynomial of the moduli space of genus zero stable maps into Grassmannians.

\subsection{Holomorphic Lefschetz formula} We will use holomorphic Lefschetz formula to determine the weights of the tangent space to a fixed point. So we recall the holomorphic Lefschetz formula (for details, see \cite{atiyah1968index}):

Let $M$ be a $G-$manifold, and $E$ be a holomorphic $G-$vector bundle on $M$. For any $g\in G$, let $M^g$ be the fixed locus of $g$, and $N^g$ be the normal bundle. Set
\begin{equation}
  \chi_g(M,E)=\sum_{q}(-1)^q\Tr_g H^q(M,E)
\end{equation}
then:
\begin{equation}
  \chi_g(M,E)=\int_{M^g}\frac{\ch_g (E|_{M^g}) \cdot \Tod (M^g)}{\ch_g \Lambda_{-1}(N^g)^*}
\end{equation}
where $\Lambda_t (E)=1+t \Lambda E+t^2\Lambda^2 E+\cdots$ for any vector bundle $E$, and $\ch_g: K_G (X)\to H^*(X,\mC)$ is the homomorphism defined in \cite{atiyah1968index}.

\section{Torus action on the moduli space}
  \subsection{Notations}
  \cite{Brion2005Lectures}Let $\vec{r}=(r_1,\cdots,r_{l+1})$ be an $(l+1)$-tuple positive integrals with $\sum_i^{l+1} r_i=k$. Let $Fl=Fl(r_1,\cdots,r_{l+1};k)=Fl(r_1,\cdots,r_{l+1})$ the flag manifold, which is the moduli space of flags of vector subspaces in $\mC^k$:
  $$
  V_1 \subset V_2 \subset \cdots \subset V_l \subset \mC^k
  $$
  with $\dim V_j=\sum_{i=1}^j r_i$. We have canonical embedding:
  \begin{equation}
    Fl(r_1,\cdots,r_{l+1})\hookrightarrow \Gr(s_1,k)\times \cdots \times \Gr(s_l,k)\hookrightarrow \mP^{N_1}\times \cdots \times \mP^{N_l}
  \end{equation}
  where $s_i=r_1+\cdots+r_i, i=1,\cdots,l$, and the second arrow is the product of Pl\"ucker embeddings. The pull backs of the hyperplane classes $H_i,i=1,\cdots,l$ form a basis of $H^2(Fl,\mZ)$ and they span the K\"{a}hler cone of $Fl(r_1,\cdots,r_{l+1})$. Let $\check{K}$ be the classes in $H_2(Fl,\mZ)$ that lies in the dual of the closure of the K\"ahler cone, and let $\{\check{H}_i|i=1,\cdots,l\}$ be the dual basis of $\{H_i|i=1,\cdots,l\}$. Then we may write $\bd\in \check{K}$ as $\bd=d_1 \check{H}_1+\cdots+d_l \check{H}_l$ with $d_i$ nonnegative integers. Let $\cM_{0}(Fl,\bd)$ denote the moduli space of genus zero stable maps into $Fl(r_1,\cdots,r_{l+1})$ of degree $\bd\in \check{K}$.

%  Over $Fl$ we have a universal flag:
%  \begin{equation}
%    D_1\subset \cdots \subset D_l \subset \underline{\mC}^k=:D_{l+1}
%  \end{equation}
%  where $D_i$ is the universal bundle of rank $\sum_{j=1}^i r_j$ and $\underline{\mC}^k$ is the trivial vector bundle. Let $c^i=1+c^i_1+\cdots c^i_{r_i}$  be the total Chern class of the dual of the quotient bundle $(D_i\slash D_{i-1})^{\vee}$, and $c=1+c_1+\cdots+c_r=\prod_{i=1}^l c_i$ be the total Chern class of $D_l$.

  \subsection{Torus action on flag manifolds}
  Let $T=\mC^*$ be an algebraic torus. Let $e_i=(0,\cdots,0,1,0,\cdots,0),i=1,\cdots,k$ be the canonical basis of $\mC^k$. $T$ acts on $\mC^k$ by:
  %, with Lie algebra $\mathfrak{t}$ canonically isomorphic to $\mR^k$. Let $e_i=(0,\cdots,0,1,0,\cdots,0),i=1,\cdots,k$ be the canonical basis of $\mathfrak{t}$. Let $\alpha_1=2\pi\sqrt{-1}\tilde \alpha_1,\cdots,\alpha_k=2\pi\sqrt{-1}\tilde \alpha_k\in \mathfrak{t}^*$ be $k$  weights of $T^k$, which defines $k$ one-dimensional representations:
%  \begin{align}
%    \rho_{i}: T^k\times \mC &\to \mC \nonumber \\
%    (t,z)&\to t^{\alpha_i}z
%  \end{align}
%  where $t^{\alpha}=t_1^{\alpha(e_1)}\cdots t_k^{\alpha(e_k)}$ for any $\alpha \in \mathfrak{t}^*$. Via the isomorphism $\mathfrak{t}\cong \mR^k$, we may identify each weight $\alpha$  with $k$-tuple of integers $(\alpha(e_1),\cdots,\alpha(e_k))\in \mZ^k$. Fix a Cartesian basis $\{e_i\}_{i=1,\cdots,k}$ of $\mC^k$, the representations $\rho_i$ induces an action on $Fl$ by acting on the total space $\mC^k$:
  \begin{align}
    T \times \mC^k &\to \mC^k \nonumber \\
    (t,(x_1,\cdots,x_k))& \mapsto (t^{\alpha_1}x_1,\cdots,t^{\alpha_k}x_k)
  \end{align}
  where $\alpha_i$ are generic integers such that $\alpha_1<\alpha_2<\cdots<\alpha_k$. This action induces a torus action on the flag manifold. For convenience, we use the matrix representation of the flag manifold: Let $M_{r,k}^{\circ}$ be the set of $r\times k$ complex matrices that has $r$ linearly independent rows. For any $ ^t(V_1,\cdots,V_r) \in M_{r,k}^{\circ}, V_i\in \mC^k$, we associate the flag $$\Span   \{V_1,\cdots,V_{r_1}\}\subset \Span\{V_1,\cdots,V_{r_1+r_2}\}\subset \cdots \subset \Span\{V_1,\cdots,V_r\}\subset \mC^k$$ in $Fl$. Thus we have a surjection (in fact a principle bundle):
  \begin{equation}
    M_{r,k}^{\circ}\to Fl(r_1,\cdots,r_{l+1})
  \end{equation}
  and $Fl(r_1,\cdots,r_{l+1})$ may be obtained as the quotient $Gl(r_1,\cdots,r_l)\backslash M_{r,k}^{\circ}$, where $Gl(r_1,\cdots,r_l)$ denotes the subgroup of $Gl(r,\mC)$ that consists of block lower triangular matrices, i.e. invertible matrices of the following form:
  \begin{equation}
  \left(
    \begin{array}{cccc}
    A_{r_1} & 0 & \cdots & 0 \\
    \ast & A_{r_2} & \cdots & 0 \\
    \vdots & \vdots &\ddots & \vdots \\
    \ast & \ast & \cdots & A_{r_l}
    \end{array}
    \right)
  \end{equation}
  $Gl(r_1,\cdots,r_l)$ acts on $M_{r,k}^{\circ}$ via left multiplication and $T$ acts in the following manner:
     \begin{equation}
  \left(
    \begin{array}{cccc}
    a_{11} & a_{12} & \cdots & a_{1k} \\
    a_{21} & a_{22} & \cdots & a_{2k} \\
    \vdots & \vdots &\vdots & \vdots \\
    a_{r1} & a_{r2} & \cdots & a_{rk}
    \end{array}
    \right)
    \cdot t :=
     \left(
    \begin{array}{cccc}
    t^{\alpha_1}a_{11} & t^{\alpha_2}a_{12} & \cdots & t^{\alpha_k}a_{1k} \\
    t^{\alpha_1}a_{21} & t^{\alpha_2}a_{22} & \cdots & t^{\alpha_k}a_{2k} \\
    \vdots & \vdots &\vdots & \vdots \\
    t^{\alpha_1}a_{r1} & t^{\alpha_2}a_{r2} & \cdots & t^{\alpha_k}a_{rk}
    \end{array}
    \right)
  \end{equation}

   We easily see from this matrix description that the fixed points of the torus action are $\{P_{I_1,\cdots,I_{l+1}}|(I_1,\cdots,I_{l+1})\in \fI(r_1,\cdots,r_{l+1})\}$, where the index set $\fI(r_1,\cdots,r_{l+1})=\{(I_1,\cdots,I_{l+1})|I_a \subset \{1,\cdots,k\},|I_a|=r_a, I_a \cap I_b=\varnothing, \forall a\neq b\}$ and $P_{I_1,\cdots,I_{l+1}}$ represents the flag:
    \begin{equation}
      \Span\{e_i|i\in I_1\} \subset \Span \{e_i|i\in I_1\cup I_2\}\subset \cdots \subset \Span \{e_i|i\in I_1\cup\cdots\cup I_l\}\subset \mC^k
    \end{equation}
    since $I_{l+1}$ is determined by $(I_1,\cdots,,I_l)$, we also denote it by $P_{I_1,\cdots,I_l;k}$.
%    corresponding to the matrix ${}^t(e_{I_1},\cdots,e_{I_l})$, in which we use the notation that for a multi-index $I=\{a_1,\cdots,a_s\}\subset \{1,\cdots,k\}$, $e_I$ represents the matrix $(e_{a_1},\cdots,e_{a_s})$.

  $Fl$ is covered by $|\fI(r_1,\cdots,r_{l+1})|=\frac{k!}{r_1!\cdots r_{l+1}!}$ affine open subsets $\{U_{I_1,\cdots,I_{l+1}}|(I_1,\cdots,I_{l+1})\in \fI(r_1,\cdots,r_{l+1})\}$, in terms of matrix, $A=(a_{ij})\in M_{r,k}^{\circ}$ lies in $U_{I_1,\cdots,I_{l+1}}$ if and only if $A^{I_1,\cdots,I_{j}}_{1,\cdots,r_1+\cdots+r_j}\neq 0, j=1,\cdots,l$, where $A_{i_1,\cdots,i_s}^{j_1,\cdots,j_s}$ is the minor with column indices $\{j_1,\cdots,j_s\}$ and row indices $\{i_1,\cdots,i_s\}$:
  \begin{equation}
    A_{i_1,\cdots,i_s}^{j_1,\cdots,j_s}=\det A\binom{j_1\cdots j_s}{i_1\cdots i_s}
  \end{equation}
  %where
%  \begin{equation}
%  A\binom{j_1\cdots j_s}{i_1\cdots i_s}=\left(
%    \begin{array}{cccc}
%    a_{i_1,j_1} & a_{i_1,j_2} & \cdots & a_{i_1,j_s}\\
%    a_{i_2,j_1} & a_{i_2,j_2} & \cdots & a_{i_2,j_s}\\
%    \vdots & \vdots & \ddots & \vdots \\
%    a_{i_s,j_1} &a_{i_s,j_2} & \cdots & a_{i_s,j_s}
%    \end{array}
%      \right)
%  \end{equation}
  The matrix representation of a flag in analogy to projective space can be viewed as  homogeneous coordinates of the flag manifold. To obtain the inhomogeneous coordinates over the affine open subset $U_{I_1,\cdots,I_{l+1}}$, observe that for any $A=(a_{ij})\in U_{I_1,\cdots,I_{l+1}}$, by linear algebra, there exists a unique $g\in Gl(r_1,\cdots,r_l)$ such that,
  \begin{equation}
    A=gB
  \end{equation}
  where $B$ is the matrix such that the submatrix $B\binom{I_1,\cdots,I_l}{1,\cdots,r}$ is an upper-diagonal matrix with diagonal elements the identities, the matrix $B$ provides the inhomogeneous coordinates. Now let us see how the torus acts on the inhomogeneous coordinates:
  \begin{align}
    A\cdot \left(
    \begin{array}{cccc}
      t^{\alpha_1}& 0 & \cdots & 0 \\
      0 & t^{\alpha_2} & \cdots & 0 \\
      \vdots & \vdots &\ddots & \vdots \\
      0 & 0 & \cdots & t^{\alpha_k}
    \end{array}
    \right) &=g\cdot  \left(
    \begin{array}{cccc}
      t^{\alpha_{I_1}} & 0 & \cdots & 0\\
      0 & t^{\alpha_{I_2}} & \cdots & 0 \\
      \vdots & \vdots & \ddots & \vdots \\
      0 & 0 & \cdots & t^{\alpha_{I_l}}
    \end{array}
    \right)\cdot \tilde{B}
  \end{align}
  where for a multi-index $I=\{i_1,\cdots,i_s\}$, $t^{\alpha_{I}}=\textrm{diag}(t^{\alpha_{i_1}},\cdots,t^{\alpha_{i_s}})$, and
  \begin{equation}
    \tilde{B}=\left(
    \begin{array}{cccc}
      t^{\alpha_{I_1}} & 0 & \cdots & 0\\
      0 & t^{\alpha_{I_2}} & \cdots & 0 \\
      \vdots & \vdots & \ddots & \vdots \\
      0 & 0 & \cdots & t^{\alpha_{I_l}}
    \end{array}
    \right)^{-1}\cdot B \cdot\left(
    \begin{array}{cccc}
      t^{\alpha_1}& 0 & \cdots & 0 \\
      0 & t^{\alpha_2} & \cdots & 0 \\
      \vdots & \vdots &\ddots & \vdots \\
      0 & 0 & \cdots & t^{\alpha_k}
    \end{array}
    \right)
  \end{equation}
  Note that the submatrix $\tilde B\binom{I_1,\cdots,I_l}{1,\cdots,r}$ is still an upper-diagonal matrix with diagonal elements the identities, so $\tilde B$ provides the inhomogeneous coordinate of $A\cdot t$. It follows that the tangent space $TFl|_{P_{I_1,\cdots,I_{l+1}}}$ at $P_{I_1,\cdots,I_{l+1}}$ as a $T$-representation splits into $\dim Fl$ one-dimensional irreducible representations with weights $\{\alpha_i-\alpha_j|i\notin I_1\cup \cdots \cup I_s,j\in I_s,s=1,\cdots,l \}$.

  %Next, we are going to lift the $T^k$ action on $Fl$ to its universal bundles $D_i,i=1,\cdots,l$. Let $V_i=\mC^{r_1+\cdots+r_i}$ be the vector space of row vectors over $\mC$, and $Gl(r_1+\cdots+r_i,\mC)$ acts on $V_i$ by right multiplication. Through the composition $Gl(r_1,\cdots,r_l)\to Gl(r_1,\cdots,r_i)\to Gl(r_1+\cdots+r_i)$, $Gl(r_1,\cdots,r_l)$ acts on $V_i$ on the right. The total space of $D_i$ is $V_i\times_{Gl(r_1,\cdots,r_l)}M_{r,k}^{\circ}$. And we lift the $T^k$ action to $D_i$ by acting only on the second factor $M_{r,k}^{\circ}$. From this description it is easy to see that over a fixed point $P_{I_1,\cdots,I_l}$, $(I_1,\cdots,I_r)\in \fI(r_1,\cdots,r_{l+1})$, $T^k$ acts on the fibre of $D_i$ with weights $\alpha_i,i\in I_1\cup\cdots\cup I_i$. And $T^k$ acts on $(D_i\slash D_{i-1})^{\vee}|_{P_{I_1,\cdots,I_l}}$ with weights $-\alpha_i, i\in I_i$.

  In summary, we have the following lemma:
  \begin{lemma}\label{lemma: weights of tangent bundle}
    If we consider $TFl|_{P_{I_1,\cdots,I_{l+1}}}$ as a $T-$representation, the weights are $\{\alpha_i-\alpha_j|i\notin I_1\cup \cdots \cup I_s,j\in I_s,s=1,\cdots,l \}$.
  \end{lemma}
\subsection{Torus action on the moduli space}
The above torus action induces a natural torus action on the moduli space of stable maps $\cM_0 (Fl,\bd)$ by acting on the target space. To see its fixed points in the moduli space. We need to know the fixed lines in $Fl(r_1,\cdots,r_{l+1})$. Since $\chi(\mP^1)=2$, every fixed line in $Fl$ connects two fixed point in $Fl$. Via the $T-$equivariant embedding:
\begin{equation}\label{eqn: natural embedding}
  Fl(r_1,\cdots,r_{l+1};k)\hookrightarrow \Gr(s_1,k)\times \cdots \times \Gr(s_l,k)
\end{equation}
Every fixed line in $Fl$ is embedded as a fixed line in the product of Grassmannians. Recall the fact that two fixed points $P_{I;k},P_{J;k}$ of $\Gr(r,k)$ is connected by a fixed line if and only if $|I\cap J|=r-1$, i.e. the index set $I$ differs from $J$ in only one element (see for example \cite{agrawal2011euler}). Explicitly, when $I=\{e_{i_1},\cdots,e_{i_{r-1}},e_a\}$ and $J=\{e_{i_1},\cdots,e_{i_{r-1}},e_b\}$, let $P_{I;k}=\Span\{e_{i_1},\cdots,e_{i_{r-1}},e_a\}$ and $P_{J;k}=\Span \{e_{i_1},\cdots,e_{i_{r-1}},e_b\}$ be the corresponding two fixed point in $\Gr(r,k)$, then up to change of coordinate, the fixed line passing through $P_{I;k}$ and $P_{J;k}$ is:
\begin{align}
  \mP^1 & \hookrightarrow \Gr(r;k) \nonumber \\
  [z:w]&\mapsto \Span \{e_{i_1},\cdots,e_{i_{r-1}},ze_a+we_b\}
\end{align}
Now we analyze the fixed lines in $Fl$, let $P_{I_1,\cdots,I_l;k}$ and $P_{\tilde I_1,\cdots,\tilde I_l;k}$ be two fixed points in $Fl(r_1,\cdots,r_{l+1})$ that is connected by a fixed line $\mP^1\hookrightarrow Fl$. Under the embedding \eqref{eqn: natural embedding}:
\begin{align}
  P_{I_1,\cdots,I_l;k}&\mapsto (P_{I_1;k},\cdots,P_{I_1\cup\cdots\cup I_l;k})\in \Gr(s_1,k)\times \cdots \times \Gr(s_l,k) \\
  P_{\tilde I_1,\cdots,\tilde I_l;k}&\mapsto (P_{\tilde I_1;k},\cdots,P_{\tilde I_1\cup\cdots\cup \tilde I_l;k})\in \Gr(s_1,k)\times \cdots \times \Gr(s_l,k)
\end{align}
and the fixed line is embedded in the product of Grassmannian. When projected to each component $\Gr(s_i,k)$, $i=1,\cdots,l$, the line becomes either a fixed line connecting $P_{I_1\cup\cdots\cup I_s;k}$ and $P_{\tilde I_1\cup\cdots\cup \tilde I_s;k}$ or a single fixed point.

Let $m$ be the smallest integer such that $P_{I_1\cup \cdots \cup I_i;k}=P_{\tilde I_1\cup\cdots \cup \tilde I_i;k}$ for $i\leq m$, and $P_{I_1\cup \cdots \cup I_{m+1};k}$ is connected to $P_{\tilde I_1\cup \cdots \cup \tilde I_{m+1};k}$ through a fixed line in $\Gr(s_{m+1},k)$. Then $I_i=\tilde I_i$ for $i\leq m$, and $I_{m+1}$ differs from $\tilde I_{m+1}$ in only one element, say $a\in I_{m+1}\setminus \tilde I_{m+1}$ and $b\in \tilde I_{m+1} \setminus I_{m+1}$.
 For convenience, we denote $I_{m+1}$ by $A\cup \{a\}$, and $\tilde I_{m+1}$ by $A\cup \{b\}$. As for the $(m+2)-$th component, since $I_1\cup\cdots \cup I_{m+2}$ must be different from $\tilde I_1\cup \cdots \cup \tilde I_{m+2}$ in at most one element, there are two possibilities: $I_{m+2}=J_{m+2}$ or $A\cup \{a\} \cup I_{m+2}=A\cup \{b\}\cup \tilde I_{m+2}$. Let $n$ be the smallest integer such that $I_{m+i+1}=\tilde I_{m+i+1}$ for $1\leq i \leq n$ and $I_{m+1}\cup \cdots \cup I_{m+2+n}=\tilde I_{m+1}\cup \cdots \cup \tilde I_{m+2+n}$. For convenience, we denote $I_{m+1+i}$ by $J_i$ , $\tilde I_{m+1+i}$ by $\tilde J_i$, $i=1,\cdots,n$, and write $I_{m+2+n}$ as $B\cup\{b\}$, and write $\tilde I_{m+2+n}$ as $B\cup \{a\}$. When $i\geq m+3+n$, $I_i$ must be the same as $\tilde I_i$, since otherwise, by a simple argument, the line will not be fixed by the torus action. Again we denote $I_{m+2+n+i}$ by $K_{i}$ and $\tilde I_{m+2+n+i}$ by $\tilde K_{i}$, $i=1,\cdots, p$, where $n+m+p+2=k$. So we may index the set of fixed lines in $Fl$ by the set
 \begin{align*}
   \cA&=\bigcup_{m+n+p+2=l+1}\cA_{m,n,p}
 \end{align*}
 where
 \begin{align*}
   \cA_{m,n,p}&=\fI_{r_1,\cdots,r_m,r_{m+1}-1,2,r_{m+2},\cdots,r_{m+n+1},r_{m+n+2}-1,r_{m+n+3},\cdots,r_{l+1}}\\
   &=\{(I_1,\cdots,I_m,A,\{a,b\},J_1,\cdots,J_n,B, K_1,\cdots,K_p)|\\ &|I_i|=r_i,|A|=r_{m+1}-1, |J_i|=r_{m+1+i},|B|=r_{m+2+n}-1,|K_i|=r_{m+n+2+i},\\
   &\textrm{ the subsets form a partition of the set} ~\{1,2,\cdots,k\}\}.
 \end{align*}
 Explicitly, for any $(I_1,\cdots,I_m,A,\{a,b\},J_1,\cdots,J_n,B, K_1,\cdots,K_p)\in \cA$, the fixed line is given by:
 \begin{align}
   [z:w]\mapsto & \Span \{e_i|i\in I_1\} \subset \cdots \Span \{e_i|i \in \cup_{1\leq j\leq m}I_j\}\nonumber \\
   &\subset \Span \{e_i,ze_a+we_b|i\in \cup_{1\leq j\leq m}I_j\cup A\}\subset \Span \{e_i,ze_a+we_b|i\in \cup_{1\leq j\leq m}I_j\cup A\cup J_1\}\subset \cdots \nonumber \\
   &\subset \Span \{e_i,ze_a+we_b|i\in \cup_{1\leq j\leq m}I_j\cup A\cup_{1\leq j\leq n}J_j\} \nonumber\\
   &\subset \Span \{e_i|i\in  \cup_{1\leq j\leq m}I_j \cup A\cup\{a,b\}\cup_{1\leq j\leq n}J_j \cup B \}\subset \cdots \nonumber \\
   &\subset \Span \{e_i|i\in  \cup_{1\leq j\leq m}I_j \cup A\cup\{a,b\}\cup_{1\leq j\leq n}J_j \cup B \cup_{1\leq j\leq p}K_j\}
 \end{align}
 which connects the two fixed point $P_{I_1,\cdots,I_m,A\cup \{a\},J_1,\cdots,J_n,\{b\}\cup B, K_1,\cdots,K_p}$ and $P_{I_1,\cdots,I_m,A\cup \{b\},J_1,\cdots,J_n,\{a\}\cup B, K_1,\cdots,K_p}$. It is obvious that the fixed line associated to every element in $\cA_{m,n,p}$ has degree $\bd=(0,\cdots,0,1,\cdots,1,0,\cdots,0)$, where the first $m$ terms and the last $p+1$ terms are zero.
 
 We remark that \cite{edwards2013genus} also contains an analysis of fixed lines in flag manifold by a slightly different argument.

\section{Computation of Poincar\'e polynomials}

\subsection{Notations}
In this section, we always assume that the weights $\alpha_1,\cdots,\alpha_{l+1}$ satisfy: $\alpha_i\gg\sum_{j<i}\alpha_j$. Let $\cA_n=\{1,2,\cdots,n\}$, and let $\cA_{n,j}$ be the collection of all subsets of $\cA_n$ with $j$ elements, $1\leq j\leq n$, For any $(I_1,\cdots,I_{l+1})\in \fI_{r_1,\cdots,r_{l+1}}$, we define the representation $V_{I_1,\cdots,I_{l+1}}$ as the direct sum of $V_{\alpha_{\mu}-\alpha_{\nu}}, \mu\notin I_1\cup\cdots\cup I_s,\nu\in I_s,s=1,\cdots, l$, where $V_{\alpha}$ is the one-dimensional representation with weight $\alpha\in \mZ$. By Lemma \ref{lemma: weights of tangent bundle} we know that $TFl|_{P_{I_1,\cdots,I_{l+1}}}=V_{I_1,\cdots,I_{l+1}}$. We also define the number $$N_{I_1,\cdots,I_{l+1}}=N_{I_1,\cdots,I_l;k}=\sum_{1\leq s \leq l}\sharp\{(i,j)|i\notin I_1\cup\cdots \cup I_s  j\in I_s,i>j\}$$, which is the number of positive weights in $V_{I_1,\cdots,I_{l+1}}$.

Note that for any $j_1,j_2,\cdots,j_{s+1}\in \mZ_+$ such that $j_1+j_2+\cdots+j_{s+1}=l+1$, we have a map $\fI_{r_1,\cdots,r_{l+1}}\to \fI_{r_1+\cdots+r_{j_1},r_{j_1+1}+\cdots+r_{j_1+j_2},\cdots,r_{j_s+1}+\cdots+r_{j_s+j_{s+1}}}$ sending $(I_1,\cdots,I_{l+1})$ to $(\cup_{1\leq i \leq j_1}I_{i},\cdots,\cup_{j_s+1\leq i \leq j_{s+1}}I_i)$. This is a fibration with fiber $\fI_{r_1,\cdots,r_{j_1}}\times \cdots \times \fI_{r_{j_s+1},\cdots,r_{j_s+j_{s+1}}}$. Hence, we have isomorphism between sets:
\begin{equation}\label{eqn: fibration}
  \fI_{r_1,\cdots,r_{l+1}}\cong\fI_{r_1+\cdots+r_{j_1},r_{j_1+1}+\cdots+r_{j_1+j_2},\cdots,r_{j_s+1}+\cdots+r_{j_s+j_{s+1}}}\times (\fI_{r_1,\cdots,r_{j_1}}\times \cdots \times \fI_{r_{j_s+1},\cdots,r_{j_s+j_{s+1}}})
\end{equation}
Let $((J_1,\cdots,J_{s+1}),(J_{1,1},\cdots,J_{1,j_1}),\cdots,(J_{s,1},\cdots,J_{s,j_{s+1}}))$ be the element belonging to the right hand side corresponding to $(I_1,\cdots,I_{l+1})$ under this isomorphism, one can easily see that:
\begin{align}
  N_{I_1,\cdots,I_{l+1}}=N_{J_1,\cdots,J_{s+1}}+N_{J_{1,1},\cdots,J_{1,j_1}}+\cdots+N_{J_{s,1},\cdots,J_{s,j_{s+1}}}
\end{align}

This fibration is very useful in our computation, as an example, we calculate $$f_{r_1,\cdots,r_{l+1}}(t)=f_{r_1,\cdots,r_l;k}(t)=\sum_{(I_1,\cdots,I_{l+1})\in \fI_{r_1,\cdots,r_{l+1}}}t^{N_{I_1,\cdots,I_{l+1}}}.$$
in fact, if we take a specific fibration  $\fI_{r_1,\cdots,r_{l+1}} \to \fI_{r_1+\cdots+r_l,r_{l+1}}$, using the corresponding isomorphism \eqref{eqn: fibration}, we have:
\begin{align}
  f_{r_1,\cdots,r_{l+1}}(t)&=f_{r_1+\cdots+r_l,r_{l+1}}(t)f_{r_1,\cdots,r_l}(t)
\end{align}
using this equation inductively, we finally have:
\begin{equation}\label{eqn: factorization of f}
  f_{r_1,\cdots,r_{l+1}}(t)=f_{r_1+\cdots+r_l,r_{l+1}}(t)f_{r_1+\cdots+r_{l-1},r_l}(t)\cdots f_{r_1,r_2}(t)
\end{equation}

\subsection{$q$-binomials}
We recall the concept of $q$-binomials, we refer the reader to the book \cite{kac2002quantum} for an beautiful exposition. For any positive integer $n$, the $q$-number of $n$ is denoted by $[n]_q:=\frac{q^n-1}{q-1}$; the $q$-factorial is defined by $[n]_q!:=[n]_q [n-1]_q \cdots [1]_q$; and the $q$-binomial $\binom{n}{k}_q:=\frac{[n]_q !}{[n-k]_q ! [k]_q !}$. We include here some basic identities:
\begin{align}
  \binom{n}{k}_q&=\binom{n}{n-k}_q \\
  \binom{n}{k}_q&=\binom{n-1}{k}_q+q^{n-k}\binom{n-1}{k-1}_q \\
  \binom{n}{k}_q&=q^k \binom{n-1}{k}_q+\binom{n-1}{k-1}_q \\
  \sum_{j=0}^a q^j \binom{d+j}{j}_q&=\binom{d+a+1}{a}_q
\end{align}
The following identity can be found in the appendix of \cite{Mart2013Poincar}:
\begin{align}
  \sum_{i+j=u}t^{i(j'+1)}\binom{i+i'}{i}_t\binom{j+j'}{j}_t=\binom{i'+j'+u+1}{u}_t \label{eqn: Identity I}
\end{align}
Using this, we have:
\begin{align}\label{eqn: identity II}
  \sum_{i+j=u}t^{i(j'+2)}\binom{i+i'}{i}_t\binom{j+j'}{j}_t&=\sum_{i+j=u}t^{i((j'+1)+1)}\binom{i+i'}{i}_t(\binom{j+j'+1}{j}_t-t^{j'+1}\binom{j+j'}{j-1}_t)\nonumber \\
  &=\binom{i'+j'+u+2}{u}_t-t^{j'+1} \sum_{i+j-1=u-1}t^{i((j'+1)+1)}\binom{i+i'}{i}_t \binom{j-1+(j'+1)}{j-1}_t \nonumber \\
  &=\binom{i'+j'+u+2}{u}_t-t^{j'+1}\binom{i'+j'+u+1}{u-1}_t
\end{align}
We have the following combinatoric interpretation of $q$-binomials:
\begin{theorem}[\cite{kac2002quantum}]\label{thm: comibinatorial discription of q-binomial}
  \begin{equation}
    \binom{n}{j}_q=\sum_{S\in \cA_{n,j}}q^{\omega(S)-j(j+1)/2}
  \end{equation}
  where $\omega (S)=\sum_{s\in S} s$.
\end{theorem}

\subsection{Poincar\'e polynomial of $Fl(r_1,\cdots,r_{l+1};k)$}
As a warm up, we compute the Poincar\'e polynomial of the flag manifold itself. By Lemma \ref{lemma: weights of tangent bundle}, the number of positive weights at the fixed point $P_{I_1,\cdots,I_l;k}$ is $N_{I_1,\cdots,I_l;k}$. By \eqref{eqn: poincare polynomial of scheme}, the Poincar\'e of $Fl$ is:
 \begin{equation}\label{eqn: poincare polynomial of flag}
   P_{Fl}(q)=\sum_{(I_1,\cdots,I_{l+1})\in \fI}q^{2N_{I_1,\cdots,I_{l+1}}}
 \end{equation}

The following lemma can be checked by direct computation:
\begin{lemma}
For any $S\in \cA_{k;r}$, we have
\begin{equation}
  N_{S;k}=\sum_{i=1}^l \omega({}^t S)-r(r+1)/2
\end{equation}
where $S(\in \cA_{n,j})\mapsto {}^t S\in \cA_{n,j}$ is the one-to-one map that maps $(a_1,\cdots,a_j)$ to $(n+1-a_j,\cdots,n+1-a_1)$.
\end{lemma}
%\begin{proof}
%  Let $S_p=(a_1\cdots,a_{r_p})$. We are going to prove that $\sharp\{(i,j)|i\notin I_1\cup\cdots \cup I_p  ,j\in I_p,i>j\}$ only depends on $S_p$ and equals to $\omega({}^tS_p)-r_p(r_p+1)/2$, then the lemma follows. In fact, if $I_1,\cdots,I_{p-1}$ lies strictly to the right of $I_p$, one can see that:
%  \begin{align*}
%    &\sharp\{(i,j)|i\notin I_1\cup\cdots \cup I_p  ,j\in I_p,i>j\} \\
%    =&\sharp\{(i,j)|i\notin I_p  ,j\in I_p,i>j\}- \sharp\{(i,j)|i\in  I_1\cup\cdots \cup I_{p-1}  ,j\in I_p,i>j\} \\
%    =&(k-a_1-(r_p-1))+(k-a_2-(r_p-2))+\cdots+(k-a_{r_p}-0))-(r_1+\cdots+r_{p-1})r_p\\
%    =&(k+1-r_1-\cdots-r_{p-1}-a_1)+\cdots+(k+1-r_1-\cdots-r_{p-1}-a_{r_p})-(1+\cdots+r_p) \\
%    =&\omega({}^tS_p)-r_p(r_p+1)/2
%  \end{align*}
%  Now keeping $S_p$ fixed, and change $S_1,\cdots,S_{p-1}$ step by step. Every time an element of $I_1\cup\cdots\cup I_{p-1}$ jumps to the left side of an element $\iota(a) \in I_p$ for some $a\in S_p$, $\sharp\{(i,j)|i\in  I_1\cup\cdots \cup I_{p-1}  ,j\in I_p,i>j\}$ decreases by one and $\sharp\{(i,j)|i\notin I_p  ,j\in I_p,i>j\}$ also decreases by one(since $\iota(a)$ increases by one), so $\sharp\{(i,j)|i\notin I_1\cup\cdots \cup I_p  ,j\in I_p,i>j\}$ keeps unchanged. And we have completed the proof.
%\end{proof}
Now, by the lemma and Theorem \ref{thm: comibinatorial discription of q-binomial}, combining \ref{eqn: factorization of f} and  \eqref{eqn: poincare polynomial of flag}, we can compute the Poincar\'e polynomial of $Fl(r_1,\cdots,r_{l+1})$:
\begin{align}
  P_{Fl}(q)&=f_{r_1,\cdots,r_{l};k}(q^2)\nonumber \\
  &=f_{r_1;k}(q^2)f_{r_2;k-r_1}(q^2)\cdots f_{r_l;k-r_1-\cdots-r_{l-1}}(q^2) \nonumber \\
  &=\binom{k}{r_1}_{q^2}\binom{k-r_1}{r_2}_{q^2}\cdots \binom{k-r_1-\cdots-r_{l-1}}{r_1}_{q^2} \nonumber \\
  &=\binom{k}{r_1 \cdots r_{l+1}}_{q^2} \nonumber
  \end{align}

\subsection{Poincar\'e polynomial of $\cM_0(Fl,\check H_i)$}
By Theorem \ref{thm: homology basis theorem for DM stack}, it suffices to compute the number of positive (or negative) weights of the tangent space at the fixed points of $\cM_0(Fl,\check H_i)$. By the analysis in the last section, the fixed point set in $\cM_0(Fl,\check H_i)$ consists of fixed lines $\mP^1\to Fl(r_1,\cdots,r_{l+1};k)$ that are parameterized by the index set $\cA_{i-1,0,l-i}$. For any  $(I_1,\cdots,I_{i-1},A,\{a,b\},B, K_1,\cdots,K_{l-i})\in \cA_{i-1,0,l-i}$, let $f:\mP^1 \to Fl(r_1,\cdots,r_{l+1};k)$ be the corresponding fixed line connecting $p_a=P_{I_1,\cdots,I_{i-1},A\cup \{a\},B\cup \{b\}, K_1,\cdots,K_{l-i}}$ and $p_b=P_{I_1,\cdots,I_{i-1},A\cup \{b\},B\cup \{a\}, K_1,\cdots,K_{l-i}}$.

Recall that since the flag manifold is convex, the moduli space is unobstructed, and the tangent space to any fixed point $(\Sigma,f:\Sigma \to Fl)$ can be identified with $\Def (\Sigma,f)$, which as a $T$-representation fits into the deformation long exact sequence:
\begin{displaymath}
\begin{array}{ccccccc}
  0&\to &\Aut(\Sigma)&\to &\Def(f)& \to &\Def(\Sigma,f) \\
  &\to &\Def(\Sigma)&\to & 0. &&
\end{array}
\end{displaymath}
where $\Def(f)=H^0(\Sigma,f^*TFl)$. Apply this long exact sequence to our case (see \cite{hori2003mirror} for analysis of stable maps to $\mP^n$), $\Def(\Sigma)=0$ and $\Aut (\Sigma)=V_0+V_{\alpha_a-\alpha_b}+V_{\alpha_b-\alpha_a}$, where $V_n$ is the one-dimensional representation with weight $n\in \mZ$. To compute weights of the representation $H^0(\Sigma,f^*TFl)$, note that $H^i(\Sigma,f^*TFl)=0,i\geq 1$, so $\Tr_g H^0(\Sigma,f^*TFl)=\chi_g (\mP^1,f^*TFl)$, and we may use the holomorphic Lefschetz formula to compute it:
\begin{align}
  \chi_g (\mP^1,f^*TFl)&=\frac{\Tr_g (TFl|_{p_a})}{1-g^{\alpha_a-\alpha_b}}+\frac{\Tr_g(TFl|_{p_b})}{1-g^{\alpha_b-\alpha_a}} \nonumber \\
  &=[(\frac{\sum_{1\leq h \leq i-1}\sum_{\nu \in I_h, \mu \notin I_1\cup\cdots\cup I_h}g^{\alpha_{\mu}-\alpha_{\nu}}}{1-g^{\alpha_a-\alpha_b}}+\frac{\sum_{1\leq h \leq i-1}\sum_{\nu \in I_h, \mu \notin I_1\cup\cdots\cup I_h}g^{\alpha_{\mu}-\alpha_{\nu}}}{1-g^{\alpha_b-\alpha_a}})] \nonumber \\
  &+[(\frac{\sum_{\mu \in B\cup K_1\cup \cdots \cup K_{l-i},\nu\in A}g^{\alpha_{\mu}-\alpha_{\nu}}}{1-g^{\alpha_a-\alpha_b}}+\frac{\sum_{\mu \in B\cup K_1\cup \cdots \cup K_{l-i},\nu\in A}g^{\alpha_{\mu}-\alpha_{\nu}}}{1-g^{\alpha_b-\alpha_a}})\nonumber \\
  &+(\frac{\sum_{\nu \in A}g^{\alpha_b-\alpha_{\nu}}}{1-g^{\alpha_a-\alpha_b}}+\frac{\sum_{\nu \in A}g^{\alpha_a-\alpha_{\nu}}}{1-g^{\alpha_b-\alpha_a}}) \nonumber \\
  &+(\frac{g^{\alpha_b-\alpha_{a}}}{1-g^{\alpha_a-\alpha_b}}+\frac{g^{\alpha_a-\alpha_{b}}}{1-g^{\alpha_b-\alpha_a}}) \nonumber \\
  &+(\frac{\sum_{\mu \in B\cup K_1\cup \cdots \cup K_{l-i}}g^{\alpha_{\mu}-\alpha_{a}}}{1-g^{\alpha_a-\alpha_b}}+\frac{\sum_{\mu \in B\cup K_1\cup \cdots \cup K_{l-i}}g^{\alpha_{\mu}-\alpha_{b}}}{1-g^{\alpha_b-\alpha_a}})]　\nonumber \\
  &+[(\frac{\sum_{\mu \in K_1\cup \cdots \cup K_{l-i},\nu\in B}g^{\alpha_{\mu}-\alpha_{\nu}}}{1-g^{\alpha_a-\alpha_b}}+\frac{\sum_{\mu \in K_1\cup \cdots \cup K_{l-i},\nu\in B}g^{\alpha_{\mu}-\alpha_{\nu}}}{1-g^{\alpha_b-\alpha_a}})\nonumber \\
  &+(\frac{\sum_{\mu \in K_1\cup \cdots \cup K_{l-i}}g^{\alpha_{\mu}-\alpha_{b}}}{1-g^{\alpha_a-\alpha_b}}+\frac{\sum_{\mu \in K_1\cup \cdots \cup K_{l-i}}g^{\alpha_{\mu}-\alpha_{a}}}{1-g^{\alpha_b-\alpha_a}})] \nonumber \\
  &+[(\frac{\sum_{1\leq h \leq l-i}\sum_{\nu \in K_h, \mu  \in K_{h+1}\cup\cdots\cup K_{l-i}}g^{\alpha_{\mu}-\alpha_{\nu}}}{1-g^{\alpha_a-\alpha_b}}+\frac{\sum_{1\leq h \leq l-i}\sum_{\nu \in K_h, \mu \in K_{h+1}\cup\cdots\cup K_{l-i}}g^{\alpha_{\mu}-\alpha_{\nu}}}{1-g^{\alpha_b-\alpha_a}})]
\end{align}
in which we use Lemma \ref{lemma: weights of tangent bundle}. Observe the following identities:
\begin{align}
  \frac{1}{1-z}+\frac{1}{1-z^{-1}}&=1, & \frac{1}{1-z}+\frac{z}{1-z^{-1}}&=1+z; \nonumber \\
  \frac{z^{-1}}{1-z}+\frac{z}{1-z^{-1}}&=z^{-1}+1+z, &\frac{1}{1-z}+\frac{z^{-1}}{1-z^{-1}}&=0 \nonumber
\end{align}
we may continue the computation:
\begin{align}
  \chi_g (\mP^1,f^*TFl)
  &=[(\sum_{1\leq h \leq i-1}\sum_{\nu \in I_h, \mu \notin I_1\cup\cdots\cup I_h}g^{\alpha_{\mu}-\alpha_{\nu}})]+[(\sum_{\mu \in B\cup K_1\cup \cdots \cup K_{l-i},\nu\in A}g^{\alpha_{\mu}-\alpha_{\nu}})\nonumber \\
  &+(\sum_{\nu \in A}(g^{\alpha_a-\alpha_{\nu}}+g^{\alpha_b-\alpha_{\nu}}))+(g^{\alpha_b-\alpha_{a}}+1+g^{\alpha_a-\alpha_b})+(\sum_{\mu \in B\cup K_1\cup \cdots \cup K_{l-i}}g^{\alpha_{\mu}-\alpha_{a}}+g^{\alpha_{\mu}-\alpha_{b}})]　\nonumber \\
  &+[(\sum_{\mu \in K_1\cup \cdots \cup K_{l-i},\nu\in B}g^{\alpha_{\mu}-\alpha_{\nu}})+0]+[(\sum_{1\leq h \leq l-i}\sum_{\nu \in K_h, \mu  \in K_{h+1}\cup\cdots\cup K_{l-i}}g^{\alpha_{\mu}-\alpha_{\nu}})]
\end{align}
So as a representation,
\begin{align}
  \Def(\Sigma,f)&=[(\sum_{1\leq h \leq i-1}\sum_{\nu \in I_h, \mu \notin I_1\cup\cdots\cup I_h}V_{\alpha_{\mu}-\alpha_{\nu}})]+[(\sum_{\mu \in B\cup K_1\cup \cdots \cup K_{l-i},\nu\in A}V_{\alpha_{\mu}-\alpha_{\nu}})\nonumber \\
  &+(\sum_{\nu \in A}(V_{\alpha_a-\alpha_{\nu}}+V_{\alpha_b-\alpha_{\nu}}))+(\sum_{\mu \in B\cup K_1\cup \cdots \cup K_{l-i}}V_{\alpha_{\mu}-\alpha_{a}}+V_{\alpha_{\mu}-\alpha_{b}})]　\nonumber \\
  &+(\sum_{\mu \in K_1\cup \cdots \cup K_{l-i},\nu\in B}V_{\alpha_{\mu}-\alpha_{\nu}})+(\sum_{1\leq h \leq l-i}\sum_{\nu \in K_h, \mu  \in K_{h+1}\cup\cdots\cup K_{l-i}}V_{\alpha_{\mu}-\alpha_{\nu}}) \nonumber \\
  &=V_{I_1,\cdots,I_{i-1},A,\{a,b\},B, K_1,\cdots,K_{l-i}}
\end{align}
So the Poincar\'e polynomial of $\cM_0(Fl,\check H_i)$ is:
\begin{align}
  P_{\cM_0(Fl,\check H_i)}(q)&=\sum_{(I_1,\cdots,I_{i-1},A,\{a,b\},B, K_1,\cdots,K_{l-i})\in \cA_{i-1,0,l-i}}q^{2 N_{I_1,\cdots,I_{i-1},A,\{a,b\},B, K_1,\cdots,K_{l-i}}} \\
&=\binom{k}{r_1, \cdots ,r_{i-1},r_i-1,2,r_{i+1}-1,r_{i+2},\cdots ,r_{l+1}}_{q^2}\\
  &=\frac{[r_i]_{q^2}[r_{i+1}]_{q^2}}{1+q^2}\frac{[k]_{q^2}}{[r_1]_{q^2}!\cdots [r_{l+1}]_{q^2}!}\\
&=\frac{[r_i]_{q^2}[r_{i+1}]_{q^2}}{1+q^2}P_{Fl}(q)
\end{align}

\subsection{Poincar\'e polynomial of $\cM_0(Fl,\check H_i+\check H_j), j-i>1$}
In this case the fixed point set of the torus action consists of maps of the form $\Sigma \to Fl(r_1,\cdots,r_{l+1};k)$ where $\Sigma$ is a nodal curve consists of two rational components such that each component is embedded in $Fl$ as a fixed line of homology class $H_i$ and $H_j$ respectively. These fixed maps are parameterized by the set $\cI_{i,j}$, where
\begin{align*}
   \cI_{i,j}&=\{(I_1,\cdots,I_{i-1},A,a,b,B, K_1,\cdots,K_{j-i-2},C,c,d,D,K_{j-i+1},\cdots,K_{l-i})|\\ &|I_j|=r_j,|A|=r_{i}-1,|B|=r_{i+1}-1,|C|=r_{j}-1,|D|=r_{j+1}-1,|K_j|=r_{i+j+1},,\\
   &\textrm{ the subsets form a partition of the set} ~\{1,2,\cdots,k\}\}.
 \end{align*}

For every $(I_1,\cdots,I_{i-1},A,a,b,B, K_1,\cdots,K_{j-i-2},C,c,d,D,K_{j-i+1},\cdots,K_{l-i})\in \cI_{i,j}$, the corresponding fixed map $f:\Sigma \to Fl$ can be described as follows:
\begin{itemize}
\item $f(p)=P_{I_1,\cdots,I_{i-1},A\cup \{a\},\{b\}\cup B, K_1,\cdots,K_{j-i-2},C\cup \{c\},\{d\}\cup D,K_{j-i+1},\cdots,K_{l-i}}$, where $p$ is the node of $\Sigma$.
\item one of the two components of $\Sigma$, say $\Sigma_1$, is mapped to the fixed line connecting $P_{I_1,\cdots,I_{i-1},A\cup \{a\},\{b\}\cup B, K_1,\cdots}$ and $P_{I_1,\cdots,I_{i-1},A\cup \{b\},\{a\}\cup B, K_1,\cdots}$
\item the other component of $\Sigma$, say $\Sigma_2$, is mapped to  the fixed line connecting $P_{\cdots, C\cup \{c\},\{d\}\cup D,K_{j-i+1},\cdots,K_{l-i}}$ and $P_{\cdots,C\cup \{d\},\{c\}\cup D,K_{j-i+1},\cdots,K_{l-i}}$
\end{itemize}
we again use the deformation long exact sequence, in this case $\Def (\Sigma)=V_{\alpha_b-\alpha_a+\alpha_d-\alpha_c}$ which corresponds to the one-dimensional deformation of the node, and $\Aut(\Sigma)=V_{\alpha_a-\alpha_b}+V_0+V_0+V_{\alpha_c-\alpha_d}$, to calculate $\Def(f)$, we use the normalization exact sequence, let $C=\mP^1 \sqcup \mP^1 \to \Sigma$ be the normalization of $\Sigma$, then we have:
\begin{equation}
  0\to \cO_{\Sigma}\to \cO_C \to \mC_{p} \to 0
\end{equation}
where $\mC_p$ is the skyscraper sheaf at the nodal point $p$. Tensoring with $f^*TFl$, we have the long exact sequence for cohomology:
\begin{align*}
  0\to H^0(\Sigma, f^* TFl)\to H^0(\Sigma_1,f^*TFl)\oplus H^0(\Sigma_2,f^*TFl) \to TFl|_{f(p)} \to 0
\end{align*}
we again use holomorphic Lefschetz formula to compute $H^0(\Sigma_1,f^*TFl)$ and $H^0(\Sigma_2,f^*TFl)$:
\begin{align}
  H^0(\Sigma_1,f^*TFl)&=V_{I_1,\cdots,I_{i-1},A,\{a,b\},B, K_1,\cdots,K_{j-i-2},C\cup\{c\},\{d\}\cup D,K_{j-i+1},\cdots,K_{l-i}}+V_0+V_{\alpha_a-\alpha_b}+V_{\alpha_b-\alpha_a}\nonumber \\
  &=V_{I_1,\cdots,I_{i-1},A,\{a,b\},B, K_1,\cdots,K_{j-i-2},C,\{c,d\},D,K_{j-i+1},\cdots,K_{l-i}}+V_{\alpha_d-\alpha_c}-V_{C,\{c\}}-V_{\{d\},D}\nonumber \\
  &+V_0+V_{\alpha_a-\alpha_b}+V_{\alpha_b-\alpha_a}
\end{align}
\begin{align}
  H^0(\Sigma_2,f^*TFl)&=V_{I_1,\cdots,I_{i-1},A\cup\{a\},\{b\}\cup B, K_1,\cdots,K_{j-i-2},C,\{c,d\},D,K_{j-i+1},\cdots,K_{l-i}}+V_0+V_{\alpha_c-\alpha_d}+V_{\alpha_d-\alpha_c}\nonumber \\
  &=V_{I_1,\cdots,I_{i-1},A,\{a,b\},B, K_1,\cdots,K_{j-i-2},C,\{c,d\},D,K_{j-i+1},\cdots,K_{l-i}}+V_{\alpha_b-\alpha_a}-V_{A,\{a\}}-V_{\{b\},B} \nonumber \\
  &+V_0+V_{\alpha_c-\alpha_d}+V_{\alpha_d-\alpha_c}
\end{align}
and
\begin{equation}
  TFl|_{f(p)}=V_{I_1,\cdots,I_{i-1},A\cup \{a\},\{b\}\cup B, K_1,\cdots,K_{j-i-2},C\cup\{c\},\{d\}\cup D,K_{j-i+1},\cdots,K_{l-i}}
\end{equation}
Hence
\begin{align}
  \Def(\Sigma,f)&=\Def(f)+\Def(\Sigma)-\Aut(\Sigma) \nonumber \\
  &=V_{I_1,\cdots,I_{i-1},A,\{a,b\},B, K_1,\cdots,K_{j-i-2},C,\{c,d\},D,K_{j-i+1},\cdots,K_{l-i}}+V_{\alpha_b-\alpha_a}+V_{\alpha_d-\alpha_c}+V_{\alpha_b-\alpha_a+\alpha_d-\alpha_c}
\end{align}
Note that by our assumption on the weights $\alpha_i$, $\alpha_b-\alpha_a+\alpha_d-\alpha_c>0$ if $b=\max \{a,b,c,d\}$ or $d=\max \{a,b,c,d\}$; and $\alpha_b-\alpha_a+\alpha_d-\alpha_c<0$ if $a=\max \{a,b,c,d\}$ or $c=\max \{a,b,c,d\}$. So the Poincar\'e polynomial is:
\begin{align}
  P_{\cM_0(Fl,\check H_i+\check H_j)}(q)&=(1+q^2+q^4+q^6)\sum_{I_1,\cdots,I_{i-1},A,\{a,b\},B, K_1,\cdots,K_{j-i-2},C,\{c,d\},D,K_{j-i+1},\cdots,K_{l-i}} \nonumber \\
  &q^{2 N_{I_1,\cdots,I_{i-1},A,\{a,b\},B, K_1,\cdots,K_{j-i-2},C,\{c,d\},D,K_{j-i+1},\cdots,K_{l-i}}} \nonumber\\
&=(1+q^2+q^4+q^6)\binom{k}{\cdots ,r_{i-1},r_i-1,2,r_{i+1}-1,r_{i+2},\cdots ,r_{j-1},r_{j}-1,2,r_{j+1}-1,r_{j+2},\cdots}_{q^2}\nonumber \\
  &=(1+q^2+q^4+q^6)\frac{[r_i]_{q^2}[r_{i+1}]_{q^2}[r_j]_{q^2}[r_{j+1}]_{q^2}}{(1+q^2)^2}\binom{k}{r_1, \cdots ,r_{l+1}}_{q^2}\nonumber\\
&=(1+q^4)\frac{[r_i]_{q^2}[r_{i+1}]_{q^2}[r_j]_{q^2}[r_{j+1}]_{q^2}}{1+q^2}P_{Fl}(q)
\end{align}

\subsection{Poincar\'e polynomial of $\cM_0(Fl,2\check H_i)$}
In this case the fixed point set of the torus action is still discrete. And it can be divided into five disjoint subsets which are indexed by $\cI_1,\cdots,\cI_5$ respectively, where
\begin{align*}
   \cI_{1}&=\fI_{r_1,\cdots,r_{i-1},r_i-1,2,r_{i+1}-1,r_{i+2},\cdots,r_{l+1}}
\end{align*}
\begin{align*}
  \cI_{2}&=\fI_{r_1,\cdots,r_{i-1},r_i-1,1,1,r_{i+1}-1,r_{i+2},\cdots,r_{l+1}}
\end{align*}
\begin{align*}
  \cI_3&= \fI_{r_1,\cdots,r_{i-1},r_i-1,1,2,r_{i+1}-2,r_{i+2},\cdots,r_{l+1}}
\end{align*}
\begin{align*}
  \cI_4&=\fI_{r_1,\cdots,r_{i-1},r_i-2,2,1,r_{i+1}-1,r_{i+2},\cdots,r_{l+1}}
\end{align*}
\begin{align*}
  \cI_5&=\fI_{r_1,\cdots,r_{i-1},r_i-2,2,1,1,r_{i+1}-2,r_{i+2},\cdots,r_{l+1}}
\end{align*}
For any $(I_1,\cdots,I_{i-1},A,\{a,b\},B,K_2,\cdots,K_{l-i})\in \cI_1$, it parameterizes a fixed map $f:\Sigma \to Fl$ described as:
\begin{itemize}
  \item $\Sigma$ is isomorphic to $\mP^1$;
  \item $f$ is a covering of degree two from $\mP^1$ to the fixed line connecting $P_{I_1,\cdots,I_{i-1},A\cup \{a\},\{b\}\cup B,K_2,\cdots,K_{l-i};k}$ and $P_{I_1,\cdots,I_{i-1},A\cup \{b\},\{a\}\cup B,K_2,\cdots,K_{l-i};k}$
\end{itemize}
By holomorphic Lefschetz formula, one can check that:
\begin{equation}\label{eqn: fI1}
  \Def (\Sigma,f)=V_{I_1,\cdots,I_{i-1},A,\{a,b\},B,K_2,\cdots,K_{l-1}}+\sum_{\nu\in A}V_{\frac{\alpha_a+\alpha_b}{2}-\alpha_{\nu}}+\sum_{\mu\in B}V_{\alpha_{\mu}-\frac{\alpha_a+\alpha_b}{2}}+V_{\alpha_b-\alpha_a}+V_{\alpha_a-\alpha_b}
\end{equation}
For any $(I_1,\cdots,I_{i-1},A,a,b,B,K_2,\cdots,K_{l-i})\in \cI_2$, it parameterizes a fixed map $f:\Sigma \to Fl$ described as:
\begin{itemize}
  \item $\Sigma$ nodal curve with two rational components $\Sigma_1$ and $\Sigma_2$;
  \item $f$ maps the node to $P_{I_1,\cdots,I_{i-1},A\cup \{a\},\{b\}\cup B,K_2,\cdots,K_{l-i};k}$;
  \item $f$ is a one-to-one map from $\mP^1$ to the fixed line connecting $P_{I_1,\cdots,I_{i-1},A\cup \{a\},\{b\}\cup B,K_2,\cdots,K_{l-i};k}$ and $P_{I_1,\cdots,I_{i-1},A\cup \{b\},\{a\}\cup B,K_2,\cdots,K_{l-i};k}$ when restricted $\Sigma_i,i=1,2$.
\end{itemize}
In this case,
\begin{equation}\label{eqn: fI2}
  \Def(\Sigma,f)=V_{I_1,\cdots,I_{i-1},A,a,b,B,K_2,\cdots,K_{l-i}}+V_{A,a}+V_{b,B}+V_{2\alpha_b-2\alpha_a}
\end{equation}
For any $(I_1,\cdots,I_{i-1},A,a,\{b_1,b_2\},B,K_2,\cdots,K_{l-i})\in \cI_3$, it parameterizes a fixed map $f:\Sigma \to Fl$ described as:
\begin{itemize}
  \item $\Sigma$ nodal curve with two rational components $\Sigma_1$ and $\Sigma_2$;
  \item $f$ maps the node to $P_{I_1,\cdots,I_{i-1},A\cup \{a\},\{b_1,b_2\}\cup B,K_2,\cdots,K_{l-i};k}$;
  \item $f$ is a one-to-one map from $\mP^1$ to the fixed line connecting $P_{I_1,\cdots,I_{i-1},A\cup \{a\},\{b_1,b_{2}\}\cup B,K_2,\cdots,K_{l-i};k}$ and $P_{I_1,\cdots,I_{i-1},A\cup \{b_{i}\},\{a,b_{3-i}\}\cup B,K_2,\cdots,K_{l-i};k}$ when restricted $\Sigma_i,i=1,2$.
\end{itemize}
In this case,
\begin{equation}\label{eqn: fI3}
  \Def(\Sigma,f)=V_{I_1,\cdots,I_{i-1},A,a,\{b_1,b_2\},B,K_2,\cdots,K_{l-i}}+V_{A,a}+V_{\alpha_{b_1}+\alpha_{b_2}-2\alpha_a}+V_{\alpha_{b_2}-\alpha_{b_1}}+V_{\alpha_{b_1}-\alpha_{b_2}}
\end{equation}
For any $(I_1,\cdots,I_{i-1},A,\{a_1,a_2\},b,B,K_2,\cdots,K_{l-i})\in \cI_4$, it parameterizes a fixed map $f:\Sigma \to Fl$ described as:
\begin{itemize}
  \item $\Sigma$ nodal curve with two rational components $\Sigma_1$ and $\Sigma_2$;
  \item $f$ maps the node to $P_{I_1,\cdots,I_{i-1},A\cup \{a_1,a_2\},\{b\}\cup B,K_2,\cdots,K_{l-i};k}$;
  \item $f$ is a one-to-one map from $\mP^1$ to the fixed line connecting $P_{I_1,\cdots,I_{i-1},A\cup \{a_1,a_2\},\{b\}\cup B,K_2,\cdots,K_{l-i};k}$ and $P_{I_1,\cdots,I_{i-1},A\cup \{a_{3-i},b\},\{a_i,\}\cup B,K_2,\cdots,K_{l-i};k}$ when restricted $\Sigma_i,i=1,2$.
\end{itemize}
In this case,
\begin{equation}\label{eqn: fI4}
  \Def(\Sigma,f)=V_{I_1,\cdots,I_{i-1},A,\{a_1,a_2\},b,B,K_2,\cdots,K_{l-i}}+V_{\alpha_{a_1}-\alpha_{a_2}}+V_{\alpha_{a_2}-\alpha_{a_1}}+V_{b,B}
\end{equation}
For any $(I_1,\cdots,I_{i-1},A,\{a_1,a_2\},b_1,b_2,B,K_2,\cdots,K_{l-i})\in \cI_5$, it parameterizes a fixed map $f:\Sigma \to Fl$ described as:
\begin{itemize}
  \item $\Sigma$ nodal curve with two rational components $\Sigma_1$ and $\Sigma_2$;
  \item $f$ maps the node to $P_{I_1,\cdots,I_{i-1},A\cup \{a_1,a_2\},\{b_1,b_2\}\cup B,K_2,\cdots,K_{l-i};k}$;
  \item $f$ is a one-to-one map from $\mP^1$ to the fixed line connecting $P_{I_1,\cdots,I_{i-1},A\cup \{a_1,a_2\},\{b_1,b_2\}\cup B,K_2,\cdots,K_{l-i};k}$ and $P_{I_1,\cdots,I_{i-1},A\cup \{b_i,a_{3-i}\},\{a_i,b_{3-i}\}\cup B,K_2,\cdots,K_{l-i};k}$ when restricted $\Sigma_i,i=1,2$, where we assume $a_1<a_2$.
\end{itemize}
In this case,
\begin{equation}\label{eqn: fI5}
  \Def(\Sigma,f)=V_{I_1,\cdots,I_{i-1},A,\{a_1,a_2\},\{b_1,b_2\},B,K_2,\cdots,K_{l-i}}+V_{\alpha_{a_1}-\alpha_{a_2}}+V_{\alpha_{a_2}-\alpha_{a_1}}+V_{\alpha_{b_1}-\alpha_{b_2}}+V_{\alpha_{b_2}-\alpha_{b_1}}
  +V_{\alpha_{b_1}+\alpha_{b_2}-\alpha_{a_1}-\alpha_{a_2}}
\end{equation}
To get the contributions of each case, we count the number of positive weights in the representations \eqref{eqn: fI1}\eqref{eqn: fI2}\eqref{eqn: fI3}\eqref{eqn: fI4}\eqref{eqn: fI5}, and then sum over the index sets $\cI_1,\cI_2,\cI_3,\cI_4,\cI_5$ respectively as in the formula \eqref{eqn: poincare polynomial of scheme}. Actually in this case, we do not have to do the complicated computation. In fact, we can use the fibrations $\pi_j: \cI_j\to \fI_{r_1,\cdots,r_{i-1},r_i+r_{i+1},r_{i+2},\cdots,r_{l+1}}$ to reduce to the Grassmannian case. For example, using the fibration $\pi_1$, the contribution of fixed locus indexed by $\cI_1$ is:
\begin{align}
  &\sum_{\substack{(I_1,\cdots,I_{i-1},J,K_2,\cdots,K_{l-i})\\ \in \fI_{r_1,\cdots,r_{i-1},r_i+r_{i+1},r_{i+2},\cdots,r_{l+1}}}}q^{2N_{I_1,\cdots,I_{i-1},J,K_2,\cdots,K_{l-i}}}\sum_{\substack{(A,\{a,b\},B)\\ \in \fI_{r_{i}-1,2,r_{i+1}-1}}}q^{2(N_{A,\{a,b\},B}+ \sharp \{\nu \in A|\frac{\alpha_a+\alpha_b}{2}>\alpha_{\nu}\}+\sharp \{\mu \in B|\frac{\alpha_a+\alpha_b}{2}<\alpha_{\mu}\}+1)} \nonumber \\
  =& \binom{k}{r_1,\cdots,r_{i-1},r_i+r_{i+1},r_{i+2},\cdots,r_{l+1}}_{q^2}\sum_{\substack{(A,\{a,b\},B)\\ \in \fI_{r_{i}-1,2,r_{i+1}-1}}}q^{2(N_{A,\{a,b\},B}+ \sharp \{\nu \in A|\frac{\alpha_a+\alpha_b}{2}>\alpha_{\nu}\}+\sharp \{\mu \in B|\frac{\alpha_a+\alpha_b}{2}<\alpha_{\mu}\}+1)}
\end{align}
and one can recognize that the summation in the second line is the contribution of fixed locus of type $\cI_1$ with the target space replaced by the Grassmannian $\Gr(r_i,r_i+r_{i+1})$. We finally get the result:
\begin{align}
  P_{\cM_0(Fl,2\check H_i)}(q)&=\binom{k}{r_1,\cdots,r_{i-1},r_i+r_{i+1},r_{i+2},\cdots,r_{l+1}}_{q^2}P_{\cM_0(\Gr(r_i,r_i+r_{i+1}),2 )}(q) \nonumber \\
  &=\frac{(1-t^{r_i})(1-t^{r_{i+1}})((1+t^{r_i+r_{i+1}})(1+t^3)-t(1+t)(t^{r_i}+t^{r_{i+1}}))}{(1-t)^2(1-t^2)^2}P_{Fl}(q)
\end{align}
where in the second equality, we use expression of $P_{\cM_0(\Gr(r_i,r_i+r_{i+1}),2 )}$ given in \cite[Theorem 3.1]{Mart2013Poincar}.

\subsection{Poincar\'e polynomial of $\cM_0(Fl,\check H_i+\check H_{i+1})$} In this case, there are three kinds of fixed maps, which are parameterized by the sets $\cI_i$, $\cI_i'$ and $\cA_{i-1,1,l-i-1}$ respectively, where
\begin{align*}
   \cI_{i}&=\fI_{r_1,\cdots,r_{i-1},r_i-1,1,1,r_{i+1}-2,1,1,r_{i+1},\cdots,r_{l+1}}
%   \{(I_1,\cdots,I_{i-1},A,a,b_1,B,b_2,c,C, K_2,\cdots,K_{l-i})|\\ &|I_j|=r_j,|A|=r_{i}-1,|B|=r_{i+1}-2,|C|=r_{i+2}-1,|K_j|=r_{i+j+1},,\\
%   &\textrm{ the subsets form a partition of the set} ~\{1,2,\cdots,k\}\}.
 \end{align*}
and
\begin{align*}
   \cI_{i}'&=\fI_{r_1,\cdots,r_{i-1},r_i-1,1,1,r_{i+1}-1,1,r_{i+2}-1,r_{i+3},\cdots,r_{l+1}}
%   \{(I_1,\cdots,I_{i-1},A,a,b_1,B,c,C, K_2,\cdots,K_{l-i})|\\ &|I_j|=r_j,|A|=r_{i}-1,|B|=r_{i+1}-1,|C|=r_{i+2}-1,|K_j|=r_{i+j+1},,\\
%   &\textrm{ the subsets form a partition of the set} ~\{1,2,\cdots,k\}\}.
 \end{align*}
For every $(I_1,\cdots,I_{i-1},A,a,b_1,B,b_2,c,C, K_2,\cdots,K_{l-i})\in \cI_{i}$, the corresponding fixed map $f:\Sigma \to Fl$ can be described as:
\begin{itemize}
\item $\Sigma$ is a nodal curve with two rational components $\Sigma_1$ and $\Sigma_2$.
\item $f(p)=P_{I_1,\cdots,I_{i-1},A\cup \{a\},\{b_1,b_2\}\cup B, \{c\}\cup C,K_2,\cdots,K_{l-i}}$, where $p$ is the node of $\Sigma$.
\item $f$ is a one-to-one map from $\mP^1$ to the fixed line connecting $P_{\cdots,I_{i-1},A\cup \{a\},\{b_1,b_2\}\cup B, \{c\}\cup C,K_2,\cdots}$ and $P_{\cdots,I_{i-1},A\cup \{b_1\},\{a,b_2\}\cup B, \{c\}\cup C,K_2,\cdots}$ when restricted to $\Sigma_1$.
\item $f$ is a one-to-one map from $\mP^1$ to the fixed line connecting $P_{\cdots,I_{i-1},A\cup \{a\},\{b_1,b_2\}\cup B, \{c\}\cup C,K_2,\cdots}$ and $P_{\cdots,I_{i-1},A\cup \{a\},\{b_1,c\}\cup B, \{b_2\}\cup C,K_2,\cdots}$ when restricted to $\Sigma_2$.
\end{itemize}
The computation of $\Def(\Sigma,f)$ is similar to that in the last case, one can easily check that:
\begin{align}
  \Def(\Sigma,f)&=\Def(f)+\Def(\Sigma)-\Aut(\Sigma) \nonumber \\
  &=V_{I_1,\cdots,I_{i-1},A,a,b_1,B, b_2,c,C,K_{2},\cdots,K_{l-i}}+V_{\alpha_{b_1}-\alpha_a+\alpha_c-\alpha_{b_2}}+V_{\alpha_{b_2}-\alpha_{b_1}}
\end{align}
For every $(I_1,\cdots,I_{i-1},A,a,B,b,c,C, K_2,\cdots,K_{l-i})\in \cI_{i}'$, the corresponding fixed map $f:\Sigma \to Fl$ can be described as:
\begin{itemize}
\item $\Sigma$ is a nodal curve with two rational components $\Sigma_1$ and $\Sigma_2$.
\item $f(p)=P_{I_1,\cdots,I_{i-1},A\cup \{a\},\{b\}\cup B, \{c\}\cup C,K_2,\cdots,K_{l-i}}$, where $p$ is the node of $\Sigma$.
\item $f$ is a one-to-one map from $\mP^1$ to the fixed line connecting $P_{\cdots,I_{i-1},A\cup \{a\},\{b\}\cup B, \{c\}\cup C,K_2,\cdots}$ and $P_{\cdots,I_{i-1},A\cup \{b\},\{a\}\cup B, \{c\}\cup C,K_2,\cdots}$ when restricted to $\Sigma_1$.
\item $f$ is a one-to-one map from $\mP^1$ to the fixed line connecting $P_{\cdots,I_{i-1},A\cup \{a\},\{b\}\cup B, \{c\}\cup C,K_2,\cdots}$ and $P_{\cdots,I_{i-1},A\cup \{a\},\{c\}\cup B, \{b\}\cup C,K_2,\cdots}$ when restricted to $\Sigma_2$.
\end{itemize}
One can check that:
\begin{equation}
  \Def(\Sigma,f)=V_{I_1,\cdots,I_{i-1},A,a,B,b,c,C,K_2,\cdots,K_{l-i}}+V_{b,B}+V_{\alpha_c-\alpha_a}
\end{equation}

For every $(I_1,\cdots,I_{i-1},A,\{a,b\},J,B, K_1,\cdots,K_{l-i-1})\in \cA_{i-1,1,l-i-1}$, the corresponding fixed map $f:\Sigma \to Fl$ can be described as:
\begin{itemize}
\item $\Sigma$ is isomorphic to $\mP^1$.
\item $f$ is a one-to-one map from $\mP^1$ to the fixed line connecting $P_{\cdots,I_{i-1},A\cup \{a\},J,\{b\}\cup B,K_1,\cdots,K_{l-i-1}}$ and $P_{\cdots,I_{i-1},A\cup \{b\},J,\{a\}\cup B, K_1,\cdots,K_{l-i-1}}$.
\end{itemize}
One can check that:
\begin{equation}
  \Def(\Sigma,f)=V_{I_1,\cdots,I_{i-1},A,\{a,b\},J,B,K_1,\cdots,K_{l-i-1}}+V_{J,\{a,b\}}
\end{equation}

To simplify the computation, we now assume the flag is a complete flag, i.e. $r_i=1, i=1,\cdots,l+1$. In this special case, $\cI_i=\varnothing$, and we only need to take summation over $\cI_i'$ and $\cA_{i-1,1,l-i-1}$. Using the fibration $\pi_1: \cI_i'\to \fI_{r_1,\cdots,r_i-1,r_i+r_{i+1}+r_{i+2},r_{i+3},\cdots,r_{l+1}}$, the contribution of $\cI_i'$ is:
\begin{align}
  &\binom{k}{1,\cdots,1,3,1,\cdots,1}_{q^2} \sum_{(a,b,c)\in \fI_{1,1,1}}q^{2(N_{a,b,c}+\delta(c-a))} \nonumber \\
  =&\binom{k}{1,\cdots,1,3,1,\cdots,1}_{q^2}(1 +2t+3t^3+t^4)\nonumber \\
  &=\frac{1+2t+3t^3+t^4}{(1+t)(1+t+t^2)}P_{Fl}(q)
\end{align}
Similarly the contribution of $\cA_{i-1,1,l-i-1}$ is:
\begin{align}
  \binom{k}{1,\cdots,1,3,1,\cdots,1}_{q^2}\binom{3}{2}t^2=\frac{3t^2}{(1+t)(1+t+t^2)}P_{Fl}(q)
\end{align}
In sum, we have:
\begin{equation}
  P_{\cM_0(Fl,\check H_i+\check H_{i+1})}(q)=\frac{1+2t+3t^2+3t^3+t^4}{(1+t)(1+t+t^2)}P_{Fl}(q)
\end{equation}

\bibliographystyle{plain}
%\bibliography{mybibli}

\begin{thebibliography}{10}

\bibitem{agrawal2011euler}
Shishir Agrawal.
\newblock The euler characteristic of the moduli space of stable maps into a
  grassmannian.
\newblock Undergraduate thesis, University of California, San Diego, 2011.

\bibitem{atiyah1968index}
Michael~F Atiyah and Isadore~M Singer.
\newblock {The index of elliptic operators: III}.
\newblock {\em The Annals of Mathematics}, 87(3):546--604, 1968.

\bibitem{Behrend2005Cohomology}
K.~Behrend.
\newblock Cohomology of stacks.
\newblock {\em Advances in Mathematics}, 198:583--622, 2005.

\bibitem{Bialynicki1973Some}
A.~Bialynicki-Birula.
\newblock Some theorems on actions of algebraic groups.
\newblock {\em Annals of Mathematics}, 98(3):480--497, 1973.

\bibitem{Bia1976Some}
A.~Bialynicki-Birula.
\newblock Some properties of the decompositions of algebraic varieties
  determined by actions of a torus.
\newblock {\em Bull.acad.polon.sci.sér.sci.math.astronom.phys},
  24(9):667--674, 1976.

\bibitem{Brion2005Lectures}
Michel Brion.
\newblock Lectures on the geometry of flag varieties.
\newblock In {\em Topics in cohomological studies of algebraic varieties},
  pages 33--85. Springer, 2005.

\bibitem{Carrell2002Torus}
James~B. Carrell.
\newblock {\em Torus Actions and Cohomology}.
\newblock Springer Berlin Heidelberg, 2002.

\bibitem{Chen2000Poincare}
Linda Chen.
\newblock Poincare polynomials of hyperquot schemes.
\newblock {\em Mathematische Annalen}, 319(2):págs. 235--252, 2000.

\bibitem{edwards2013genus}
Gregory Edwards.
\newblock On genus zero stable maps to the flag variety.
\newblock Undergraduate thesis, University of California, San Diego, 2013.

\bibitem{Fulton1996Notes}
W.~Fulton and R.~Pandharipande.
\newblock Notes on stable maps and quantum cohomology.
\newblock {\em Arxiv Cornell University Library}, 1996.

\bibitem{hori2003mirror}
Kentaro Hori.
\newblock {\em Mirror symmetry}, volume~1.
\newblock American Mathematical Soc., 2003.

\bibitem{kac2002quantum}
Victor Kac and Pokman Cheung.
\newblock {\em Quantum calculus}.
\newblock Springer Science \& Business Media, 2002.

\bibitem{Manin1998Stable}
Yuri~I. Manin.
\newblock Stable maps of genus zero to flag spaces.
\newblock {\em Topological Methods in Nonlinear Analysis}, (2):207--217, 1998.

\bibitem{Mart2013Poincar}
Alberto~L\'opez Mart\'in.
\newblock {Poincar\'e polynomials of stable map spaces to Grassmannians}.
\newblock In {\em Rendiconti del Seminario matematico della Universit\`{a} di
  Padova}, pages 193--208, 2013.

\bibitem{Oprea2006Tautological}
Dragos Oprea.
\newblock Tautological classes on the moduli spaces of stable maps to pr via
  torus actions.
\newblock {\em Advances in Mathematics}, 207(2):661--690, 2006.

\bibitem{Skowera2013Bia}
Jonathan Skowera.
\newblock Bialynicki-birula decomposition of deligne-mumford stacks.
\newblock {\em Proceedings of the American Mathematical Society},
  141(6):1933--1937, 2013.

\bibitem{Str1987On}
Stein~Arild Str{\o}mme.
\newblock {\em On parametrized rational curves in Grassmann varieties}.
\newblock Springer Berlin Heidelberg, 1987.

\end{thebibliography}

\end{document}